\documentclass[11pt,a4paper,english]{article}

\usepackage[T1]{fontenc}
\usepackage[utf8]{inputenc}
\usepackage{babel}
\usepackage[final]{microtype} 
\usepackage{lmodern}

\usepackage[a4paper, hmargin={2.7cm,2.7cm},vmargin={3.3cm,3.3cm}]{geometry}
\usepackage{fancyhdr, lastpage}

\usepackage{mathtools} 
\usepackage{amssymb, amsfonts, amsthm, amsxtra, bm}
\usepackage{mathrsfs}  
\usepackage{dsfont}    

\usepackage{graphicx}
\graphicspath{{Figures/}}
\usepackage{tikz}
\usetikzlibrary{calc}
\usetikzlibrary{patterns}
\usepackage{pgfplots}
\usepgfplotslibrary{fillbetween}
\pgfplotsset{compat=newest}

\definecolor{sienna}{RGB}{160,82,45}
\definecolor{violet}{RGB}{238,130,238}
\definecolor{lightorange}{RGB}{253,208,162}

\usepackage{caption}
\usepackage[margin=0.2cm]{subcaption}
\usepackage{float}

\usepackage{enumitem} 
\setlist{topsep=0.4em, partopsep=0.2em, itemsep=0.1em, parsep=0.05em}
\usepackage{url}
\usepackage{siunitx}
\usepackage{etoolbox} 
\usepackage{aliascnt} 
\usepackage{cite}

\usepackage[bf,compact,small]{titlesec}
\titlespacing*{\section}{0pt}{14pt}{4pt}
\titlespacing*{\subsection}{0pt}{8pt}{3pt}

\makeatletter
\patchcmd{\ttlh@hang}{\parindent\z@}{\parindent\z@\leavevmode}{}{}
\patchcmd{\ttlh@hang}{\noindent}{}{}{}
\makeatother

\pgfplotsset{every axis/.append style={
    axis x line=middle,
    axis y line=middle,
    axis line style={->},
    xlabel style={at={(ticklabel* cs:1)},anchor=north west},
}}

\DeclarePairedDelimiter{\itvoo}{\lparen}{\rparen}      
\DeclarePairedDelimiter{\abs}{\lvert}{\rvert}          
\DeclarePairedDelimiter{\norm}{\lVert}{\rVert}         
\DeclarePairedDelimiter{\ip}{\langle}{\rangle}         

\providecommand\given{} 
\newcommand\SetSymbol[1][]{%
    \nonscript\:\mathrel{:}\allowbreak\nonscript\:\mathopen{}}
\DeclarePairedDelimiterX\set[1]\lbrace\rbrace{%
   \renewcommand\given{\SetSymbol[\delimsize]}#1}

\numberwithin{equation}{section}
\allowdisplaybreaks[4]

\newtheorem{maintheorem}{Theorem}

\newaliascnt{lemma}{theorem}
\newtheorem{lemma}[lemma]{Lemma}
\aliascntresetthe{lemma}
\newaliascnt{proposition}{theorem}
\newtheorem{proposition}[proposition]{Proposition}
\aliascntresetthe{proposition}
\newaliascnt{corollary}{theorem}

\aliascntresetthe{corollary}
\theoremstyle{definition}

\newtheorem{example}{Example}
\theoremstyle{remark}

\newtheorem{remark}{Remark}

\DeclareMathOperator{\exponential}{e}

\DeclareMathOperator{\ord}{ord}
\DeclareMathOperator{\Wr}{Wr}
\DeclareMathOperator{\spanop}{span}

\newcommand{\R}{\mathbb{R}}
\newcommand{\C}{\mathbb{C}}
\newcommand{\Z}{\mathbb{Z}}
\newcommand{\N}{\mathbb{N}}
\newcommand{\Nzero}{\N_0}
\newcommand{\Q}{\mathbb{Q}}

\newcommand{\cF}{\mathcal{F}}

\newcommand{\myexp}[1]{\exponential^{#1}}

\newcommand{\G}{\mathcal G}

\makeatletter
\def\maketimestamp{\count255=\time
\divide\count255 by 60\relax
\edef\thetime{\the\count255:}%
\multiply\count255 by-60\relax
\advance\count255 by\time
\edef\thetime{\thetime\ifnum\count255<10 0\fi\the\count255}
\edef\thedate{\number\day-\ifcase\month\or Jan\or Feb\or Mar\or
             Apr\or May\or Jun\or Jul\or Aug\or Sep\or Oct\or
             Nov\or Dec\fi-\number\year}
\def\timstamp{\hbox to\hsize{\tt\hfil\thedate\hfil\thetime\hfil}}}
\makeatother
\maketimestamp

\fancypagestyle{plain}{%
\fancyhf{} 
\fancyfoot[C]{\small \thepage{} of \pageref{LastPage}}}

\makeatletter
\def\blfootnote{\xdef\@thefnmark{}\@footnotetext} 
\if@titlepage
  \renewenvironment{abstract}{%
      \titlepage \null\vfil
      \@beginparpenalty\@lowpenalty
      \begin{center}\bfseries \abstractname \@endparpenalty\@M \end{center}}%
     {\par\vfil\null\endtitlepage}
\else
  \renewenvironment{abstract}{%
      \small
      \list{}{\setlength{\leftmargin}{3em}\setlength{\rightmargin}{3em}}
      \item[\textbf{\abstractname:}]}
      {\endlist}
\fi
\makeatother

\usepackage{hyperref}
\hypersetup{
    colorlinks=true,
    linkcolor={red!50!black},
    citecolor={blue!50!black},
    urlcolor={blue!80!black},
    pdfview={FitH},
    pdfauthor={Jakob Lemvig},
    pdftitle={On the Gabor frame set of polynomial--Gaussian windows},
}

\usepackage[noabbrev,nameinlink,capitalize]{cleveref}
\crefname{maintheorem}{Theorem}{Theorems}
\crefname{mainconj}{Conjecture}{Conjectures}

\begin{document}

\title{On the Gabor frame set of polynomial--Gaussian windows}

\author{Jakob Lemvig\footnote{Technical University of Denmark. E-mail: \protect\url{jakle@dtu.dk}}\phantom{$\ast$}}

\date{\today}

\blfootnote{2020 {\it Mathematics Subject Classification.} Primary 42C15. Secondary: 30D15, 42C40}
\blfootnote{{\it Key words and phrases.} Gabor frame, frame set, Hermite function, polynomial--Gaussian window,
Wronskian, shift-invariant space. }

\maketitle
\thispagestyle{plain}

\begin{abstract}
We prove two sufficient conditions for the frame property of Gabor systems
$\mathcal{G}(g,a,b)=\{e^{2\pi i mb\,\cdot}g(\cdot-ak):k,m\in\mathbb{Z} \}$ generated by
polynomial--Gaussian windows $g(t)=P(t)e^{-ct^2}$, where $c>0$ and $P$ is a
nonzero complex polynomial of degree $n$. The conditions depend only on
the product $ab$ of the lattice parameters and the degree of $P$.
First, $\mathcal{G}(g,a,b)$ is a frame whenever $ab<1$ is irrational.
Second, if $ab=p/q<1$, with $p,q\in\N$ coprime, then $\mathcal{G}(g,a,b)$ is a
frame whenever $q-p>n$. Consequently, in the region $ab<1$, failure of
the frame property can occur only on rational hyperbolas whose reduced
denominator and numerator differ by at most $n$.
The proof extends the translated-Wronskian method of Faulhuber and
Petersen [\emph{The frame set of the first Hermite function},
arXiv:2609.02610] from the first Hermite function to arbitrary polynomial
multiples of a Gaussian, combining a multiplicity estimate for
Wronskians with Jensen's formula.
\end{abstract}

\section{Main results}
\label{sec:results}

For $g\in L^2(\R)$ and $a,b>0$, the Gabor system generated by $g$ is
\begin{equation*}
  \G(g,a,b)=\set{\myexp{2\pi ibm t}g(t-ak) \given k,m\in\Z }.
\end{equation*}
The Gabor system $\G(g,a,b)$ is a frame for $L^2(\R)$ if there exist $A,B>0$ such that
\begin{equation}
  A\norm{f}_2^2\le
  \sum_{k,m\in\Z}\abs{\ip{f, \myexp{2\pi ibm t}g(t-ak)}}^2
  \le B\norm{f}_2^2,\qquad f\in L^2(\R).
 \label{eq:frame}
\end{equation}
The frame set of $g$ is
\begin{equation*}
  \cF(g)=\{(a,b)\in(0,\infty)^2:\G(g,a,b)
        \text{ is a frame for }L^2(\R)\}.
\end{equation*}
The frame set $\cF(g)$ is precisely the set of pairs of translation and modulation parameters $(a,b)$ for which $\G(g,a,b)$ yields stable analysis and unconditional reconstruction of signals in $L^2(\R)$. This fact explains the importance of the frame set in Gabor analysis.

Define \emph{polynomial--Gaussian windows} to be functions of the form
\begin{equation}
  g(t)=P(t)\myexp{-ct^2},\qquad c>0,\qquad
  0\ne P\in\C[t],\qquad n=\deg P.
 \label{eq:window}
\end{equation}
We impose no restrictions on $P$ other than being a nonzero polynomial of degree $n$; in particular, we make no assumptions on the parity of $P$, the location of its zeros, or the signs or arguments of its coefficients. We will study the frame set $\cF(g)$ of polynomial--Gaussian windows.

Our main examples of polynomial--Gaussian windows are the Hermite functions, which are Hermite polynomials weighted by a Gaussian. They are defined by
\begin{equation}
  h_n(t)=\kappa_n(-1)^n \myexp{\pi t^2}
               \frac{d^n}{dt^n}\myexp{-2\pi t^2},\qquad n\in\Nzero,
  \label{eq:hermite}
\end{equation}
where $\kappa_n>0$ is chosen so that $\norm{h_n}_2=1$. Hermite functions arise in many different contexts, e.g., as eigenfunctions of several fundamental operators in analysis. Most relevant for us is that they appear 
as eigenfunctions of the harmonic oscillator
\[
  H=-\frac{d^2}{dt^2}+(2\pi t)^2,
  \qquad
  Hh_n=2\pi(2n+1)h_n.
\] 
The Gaussian function is the $0$th Hermite function $h_0(t)=2^{1/4}\myexp{-\pi t^2}$. For $c=\pi$, the polynomial--Gaussian windows of degree at most $n$ are
precisely the nonzero elements of $\spanop\{h_0,\ldots,h_n\}$, while other values of the Gaussian decay rate $c>0$ are obtained by dilation of this set of functions.

Every polynomial--Gaussian window belongs to the Feichtinger algebra
$M^1(\R)$. For any nonzero $g\in M^1(\R)$, the frame set $\cF(g)$ is an
open subset of $\{(a,b)\in(0,\infty)^2:ab<1\}$ and contains
$(0,\varepsilon)^2$ for some $\varepsilon>0$; see, e.g., \cite{LyubarskiiGabor2013}. Thus, the frame property of any polynomial--Gaussian window is only possible in the so-called oversampling regime where the time-frequency lattice product $ab$ is strictly less than one.

Our first result shows that irrationality of the lattice product $ab$ in the region $ab<1$ is sufficient for the frame property of any polynomial--Gaussian window, independently of the degree.
\begin{maintheorem}
  \label{thm:irrational}
  Let $g$ be a polynomial--Gaussian window as in \eqref{eq:window}. If $a,b>0$ satisfy
  \begin{equation}
    0<ab<1,\qquad ab\notin\Q,
    \label{eq:irrational-condition}
  \end{equation}
  then $\G(g,a,b)$ is a frame, i.e., $(a,b)\in\cF(g)$.
\end{maintheorem}

For rational products, the degree enters through the difference between
the reduced denominator and numerator.
\begin{maintheorem}
  \label{thm:rational}
  Let $g$ be a polynomial--Gaussian window as in \eqref{eq:window} with $\deg P=n\in\Nzero$. Let $a,b>0$. 
  Suppose $ab \in \itvoo{0,1} \cap \Q$. Let $p,q\in\N,\gcd(p,q)=1$ be such that $ab=p/q$.  If
  \begin{equation}
    q-p>n,
    \label{eq:degree-gap}
  \end{equation}
  then $\G(g,a,b)$ is a frame.
\end{maintheorem}

To express the two theorems together, set $\mathcal{E}_0=\emptyset$ and, for $n\in\N$,
\begin{equation*}
  \mathcal{E}_n=\set*{\frac pq\in(0,1) \given p,q\in\N,\, \gcd(p,q)=1,
                                   \, 1\le q-p\le n }.
\end{equation*}
Together, \Cref{thm:irrational,thm:rational} read
\begin{equation}
  \set*{(a,b)\in(0,\infty)^2 \given ab<1,\, ab\notin\mathcal{E}_n }
  \subset \cF(g).
  \label{eq:combined}
\end{equation}
In particular, if $\G(g,a,b)$ with $ab<1$ is not a frame, then $ab$ must
belong to $\mathcal{E}_n$, where $n$ is the degree of the polynomial part of $g$.

\subsection{Known results and related work}
\label{sec:related-work}

The study of frame sets for Hermite windows connects Gabor analysis with
sampling in Bargmann--Fock spaces, the Bargmann transform, and the theory of theta and elliptic functions. Our results unify and extend a number of important results in the literature on the frame set of Hermite functions \cite{LyubarskiuiFrames1992,SeipDensity1992a,FaulhuberPetersen2026,GroechenigGabor2009} as detailed below.

The frame set of the Gaussian was characterized by Lyubarskii and by
Seip and Wallst\'en \cite{LyubarskiuiFrames1992,SeipDensity1992a}. For
separable time-frequency lattices $a\Z\times b\Z$ it is exactly the region $\set{(a,b) \given ab<1}$. Taking $n=0$ in
\eqref{eq:combined} recovers this characterization:
\begin{equation*}
  \cF(h_0)=\set{(a,b)\in(0,\infty)^2 \given ab<1}.
\end{equation*}

Gr\"ochenig and Lyubarskii \cite{GroechenigGabor2007a,GroechenigGabor2009} proved that
\begin{equation}
  \set*{(a,b)\in(0,\infty)^2 \given ab<\frac1{n+1}} \subset \cF(g)
  \label{eq:GL-bound}
\end{equation}
for every polynomial--Gaussian window $g$ of degree at most $n$.
Indeed, \cite[Proposition~3.4]{GroechenigGabor2007a} treats
nonzero linear combinations of the first $n+1$ Hermite functions.
Thus, the statement for $n$th-order polynomial--Gaussian windows with 
arbitrary Gaussian parameter $c>0$ follows
by dilation.
Our results recover this sufficient region since, if
$ab=p/q<1/(n+1)$, then $q-p>np\ge n$, so
\Cref{thm:rational} applies, while irrational products are covered
by \Cref{thm:irrational}.

Lyubarskii and Nes \cite{LyubarskiiGabor2013} showed that an odd window in the
Feichtinger algebra $M^1(\R)$ cannot generate a frame when
\begin{equation*}
  ab=\frac{q-1}{q},\qquad q\in\N,\quad q\ge2.
\end{equation*}
This applies, in particular, to every odd Hermite function. For $h_1$,
combining this obstruction with \eqref{eq:combined} gives
\begin{equation*}
  \cF(h_1)=\set{(a,b)\in(0,\infty)^2 \given
     ab<1,\quad ab\notin\mathcal{E}_1},
  \qquad
  \mathcal{E}_1=\set{\frac{q-1}{q} \given q\in\N,\ q\ge2}.
\end{equation*}
Thus our results also recover the first-Hermite characterization of
Faulhuber and Petersen \cite{FaulhuberPetersen2026} that was conjectured in \cite{LyubarskiiGabor2013}.

Horst, Videb\ae k, and the author \cite{HorstNonframe2025,Lemvigsome2017} uncovered additional obstructions for the frame property of Hermite functions of order at least two, disproving the earlier
frame-set conjecture~\cite{Groechenigmystery2014} for every Hermite function of order at least two. In particular, for each $n\ge 2$, there exist pairs $(a,b)$ on the hyperbolas $ab=1/2$, $ab=1/3$, and $ab=2/3$, for which $\G(h_n,a,b)$ is not a frame. For odd $n\ge 3$, there are obstructions to the frame property of $\G(h_n,a,b)$ on the hyperbolas $ab=1/3$ and $ab=1/4$. These obstructions are all rational products with reduced denominator and numerator differing by at most $n$, in agreement with \Cref{thm:rational}, showing that the sufficient condition \eqref{eq:degree-gap} is sharp for $n=2$ and $n=3$ as $q-p=n$ for $ab=1/3$ and $ab=1/4$, respectively. For higher orders $n>2$, there are several obstructions to the frame property on each hyperbola, see \cite{HorstNonframe2025}.

Faulhuber, Shafkulovska, and Zlotnikov \cite{FaulhuberFrame2025} proved that $((n+1)^{-1/2},(n+1)^{-1/2}) \in \cF(h_n)$ for every
$n\ge4$. Their subsequent work \cite{FaulhuberAsymptotic2026} establishes, for every
$\eta>0$ and all sufficiently large $n$, the sufficient condition
$ab\le n^{-2/3-\eta}$. For sufficiently large $n$, \cite{FaulhuberAsymptotic2026} also establishes frame regions near the coordinate axes (i.e., when either $a$ or $b$ is sufficiently small) with products $ab$ bounded strictly below $1/2$.
These results include rational products not necessarily covered by
\eqref{eq:degree-gap} and are complementary to the degree-based
criterion in \Cref{thm:rational}.

The closest work to our proof is the paper of Faulhuber and
Petersen \cite{FaulhuberPetersen2026}. Their argument combines the pre-Gramian injectivity characterization of Gr\"ochenig, Romero, and St\"ockler \cite{GroechenigSampling2018}
with Wronskians of entire functions formed from translates of a Gaussian. 
We use the same Wronskian construction, but we add an additional step using an
order-$n$ multiplicity lemma (\Cref{sec:wronskian}), followed, in \Cref{sec:proofs},
by a direct application of Jensen's formula and a Gaussian zero-counting principle, closely related to \cite[Theorem~4.3]{GroechenigSharp2020}.

Liehr, Shafkulovska, and Taylor \cite{LiehrDerivatives2026} establish
the first-Hermite characterization for time-frequency lattices $\Lambda \subset \R^2$ with \emph{rational} determinant,
including nonseparable lattices, i.e., lattices that are not of the form $a\Z\times b\Z$ for some $a,b>0$. Corollary~5.1 in \cite{LiehrDerivatives2026} also shows that
for $n\ge1$ and $g\in\spanop\{h_0,\ldots,h_n\}\setminus\{0\}$ the rational
condition
\begin{equation}
  q\ge np+2
  \label{eq:LST-condition}
\end{equation}
is sufficient for the frame property of $\G(g,a,b)$ when $ab=p/q$ with $p,q\in\N$ coprime.
Restricted to separable lattices, \Cref{thm:rational}
replaces \eqref{eq:LST-condition} by $q\ge p+n+1$. The difference between the two
thresholds on $q$ is
\begin{equation*}
  (np+2)-(p+n+1)=(n-1)(p-1).
\end{equation*}
The thresholds agree when $n=1$ or $p=1$, while the
threshold in \Cref{thm:rational} is strictly smaller when $n,p>1$. 

Finally, Ulanovskii and Zlotnikov \cite{UlanovskiiPeriodic2026}
prove that the bound $ab<1/(n+1)$ in \eqref{eq:GL-bound} is sharp
as a sufficient product threshold on $ab$. Indeed, for each $n\ge1$ and each
separable lattice of product $ab=1/(n+1)$, they show that there exists a nonzero polynomial--Gaussian window
of degree at most $n$ that does not
generate a frame. Here $q-p=n$, which is exactly an endpoint left
open by \Cref{thm:rational}, showing the sharpness, for every $n\ge1$, of the sufficient condition \eqref{eq:degree-gap} when $p=1$.

\subsection{Outline}

In \Cref{sec:reduction} we reduce failure of the frame property to a
constant-coefficient differential relation for an entire function
formed from translates of a Gaussian. \Cref{sec:wronskian} proves the
multiplicity lemma and the nonvanishing, quasi-periodicity,
and growth of the translated Wronskian. In \Cref{sec:proofs}, a
zero-counting estimate gives both main theorems.

\section{Reduction to a differential relation}
\label{sec:reduction}

We first record a pre-Gramian criterion for separable lattices and
$M^1$ windows, in the spirit of the injectivity characterizations of
Gabor frames in \cite{GroechenigGabor2007}. It follows from
\cite[Theorems~2.1 and~3.3]{GroechenigSampling2018}, together with
Ron--Shen duality, as shown below.

\begin{proposition}[Pre-Gramian characterization]
\label{prop:pregramian}
Let $g\in M^1(\R)$ and $a,b>0$. For $x\in\R$, the pre-Gramian
\begin{equation*}
  P_g(x)
  =\bigl(g(x+aj-k/b)\bigr)_{j,k\in\Z}
\end{equation*}
defines a bounded operator on $\ell^\infty(\Z)$.
The system $\G(g,a,b)$ is a frame if and only if
$P_g(x)$ is injective on $\ell^\infty(\Z)$ for every $x\in\R$.
Equivalently, $\G(g,a,b)$ is not a frame if and only if there exist
$x\in\R$ and $\bm{0}\ne\bm{c}=(c_k)_{k\in\Z}\in\ell^\infty(\Z)$ such that
\begin{equation}
  \sum_{k\in\Z}c_k g(x+aj-k/b)=0,
  \qquad \text{for all } j\in\Z.
  \label{eq:bounded-coefficient-relation}
\end{equation}
\end{proposition}

\begin{proof}
Put $\alpha=ab$ and $h(t)=b^{-1/2}g(t/b)$. Since the dilation
$f(t)\mapsto b^{-1/2}f(t/b)$ is unitary on $L^2(\R)$, and it maps $\G(g,a,b)$ onto
$\G(h,\alpha,1)$, the Gabor system $\G(g,a,b)$ has the same frame property as $\G(h,\alpha,1)$.
Moreover,
\begin{equation}
  h(y+\alpha j-k)
  =b^{-1/2}g(y/b+aj-k/b).
  \label{eq:pregramian-scaling}
\end{equation}
Hence, it suffices to consider the Gabor system $\G(h,\alpha,1)$ generated by $h$ and the lattice
$\alpha\Z\times\Z$.

We have $h\in M^1(\R)\subset W_0(\R)$, where
\[
  W_0(\R)=\set*{f\in C(\R) \given \norm{f}_{W_0} :=
    \sum_{k\in\Z}\sup_{t\in[k,k+1]}|f(t)|<\infty }.
\]
The absolute row sums of the matrix $\bigl[h(y+\alpha j-k)\bigr]_{j,k\in\Z}$ are uniformly bounded in $y\in\R$ and $j \in\Z$ by $\norm{h}_{W_0}$, so the pre-Gramian $P_g$ defines a bounded operator on $\ell^\infty(\Z)$ by Schur's test.

Suppose first that $h$ does not have stable integer shifts, meaning
that $\{h(\cdot-k):k\in\Z\}$ is not a Riesz sequence.
By \cite[Theorem~2.1]{GroechenigSampling2018}, there is a nonzero
$\bm{c}\in\ell^\infty(\Z)$ satisfying
\[
  \sum_{k\in\Z}c_k h(t-k)=0,\qquad \text{for all }t\in\R.
\]
Thus every pre-Gramian has a nontrivial kernel. Suppose towards a contradiction that $\G(h,\alpha,1)$ is a frame. By Ron--Shen duality \cite{RonDuality1997}, the adjoint system $\G(h,1,\alpha^{-1})$ is then a Riesz sequence, and hence its subsystem of integer translates $\{h(\cdot-k):k\in\Z\}$ is a Riesz sequence, contradicting the assumption that $h$ does not have stable integer shifts. Thus, $\G(h,\alpha,1)$ is not a frame.

It remains to consider $h$ with stable integer shifts. For
$\Lambda=\alpha\Z$, the collection $W(\Lambda)$ of weak limits of
real translates of $\Lambda$ is
\[
  W(\alpha\Z)=\{y+\alpha\Z:y\in\R\}.
\]
Indeed, since $x_m+\alpha\Z=\tilde x_m+\alpha\Z$ with $\tilde x_m\in[0,\alpha]$, compactness of $[0,\alpha]$ yields a convergent subsequence $\tilde x_{m_\ell}\to y$, so any weak limit of translates of $\alpha\Z$ is of the form $y+\alpha\Z$. Conversely, every translate $y+\alpha\Z$ is a constant weak limit.
Since $-\Lambda=\Lambda$,
\cite[Theorem~3.3]{GroechenigSampling2018} now says that
$\G(h,\alpha,1)$ is a frame if and only if
$\bigl(h(y+\alpha j-k)\bigr)_{j,k\in\Z}$ is injective on
$\ell^\infty(\Z)$ for every $y\in\R$.
Equation~\eqref{eq:pregramian-scaling} relates this to injectivity of the pre-Gramian $P_g$ on $\ell^\infty(\Z)$, which completes the proof.
\end{proof}

The use of $\ell^\infty$ is crucial as injectivity on $\ell^2$ alone would not imply a positive lower frame bound. Note, however, that because this compactness-based characterization is qualitative, our approach yields no explicit estimates for the lower frame bound $A$ in \eqref{eq:frame}.

We now extend the Gaussian reduction in
\cite[Proposition~2.5]{FaulhuberPetersen2026} from the first Hermite
function to arbitrary polynomial--Gaussian windows, replacing the first
derivative by a constant-coefficient differential operator of order $n$. Its proof is based on the simple fact that the $r$th derivative of a Gaussian is a polynomial of degree $r$ times the same Gaussian. 

\begin{proposition}[Gaussian reduction]
\label{prop:gaussian-reduction}
Let $g(t)=P(t)\myexp{-ct^2}$ be a polynomial--Gaussian window
as in \eqref{eq:window} with $\deg P=n$,
and let $a,b>0$. Set
\begin{equation*}
  \alpha=ab,\qquad \beta=c/b^2,\qquad
  \phi_\beta(z)=\myexp{-\beta z^2},\qquad D=\frac{d}{dz}.
\end{equation*}
There is a unique polynomial $Q(z)=\sum_{r=0}^n q_r z^r$ of degree $n$
such that
\begin{equation}
  Q(D)\phi_\beta(t)=P(t/b)\phi_\beta(t)=g(t/b),
  \qquad t\in\R.
  \label{eq:differential-window}
\end{equation}
If $\G(g,a,b)$ is not a frame, there exist $x_*\in\R$ and
$\bm{0}\ne\bm{c}=(c_k)_{k\in\Z}\in\ell^\infty(\Z)$ such that
\begin{equation}
  F(z)=\sum_{k\in\Z}c_k \myexp{-\beta(z-k)^2},\qquad z\in\C,
  \label{eq:gaussian-entire-function}
\end{equation}
defines a nonzero entire function satisfying
\begin{equation}
  Q(D)F(x_*+j\alpha)
  =\sum_{r=0}^n q_rF^{(r)}(x_*+j\alpha)=0,
  \qquad \text{for all } j\in\Z,
  \label{eq:lattice-differential-relation}
\end{equation}
and the imaginary quasi-periodicity relation
\begin{equation}
  F\Bigl(z+\frac{i\pi}{\beta}\Bigr)
  =\myexp{\pi^2/\beta-2\pi iz}F(z),\qquad z\in\C.
  \label{eq:F-automorphy}
\end{equation}
Here the series for $F$ and all its termwise derivatives converge
absolutely and locally uniformly on $\C$. Moreover, for every
$r\in\Nzero$ there is a constant $C_{r,\beta}>0$ such that
\begin{equation}
  |F^{(r)}(t+is)|
  \le C_{r,\beta}\|\bm{c}\|_{\ell^\infty}
       (1+|s|)^r \myexp{\beta s^2},\qquad t,s\in\R.
  \label{eq:F-growth}
\end{equation}
\end{proposition}

\begin{proof}
Define the polynomials $R_r(z)$, $r\ge0$, recursively by $R_0=1$ and $R_{r+1}(z)=R_r'(z)-2\beta zR_r(z)$. By induction, $D^r\phi_\beta=R_r\phi_\beta$, and $R_r$ has degree $r$ with leading coefficient $(-2\beta)^r\ne0$.    
Hence, $R_0,\ldots,R_n$ form a
basis of the polynomials of degree at most $n$, so expanding $P(z/b)=\sum_{r=0}^n q_r R_r(z)$
in this basis gives a unique polynomial $Q(z)=\sum_{r=0}^n q_r z^r$ of degree $n$ satisfying
\eqref{eq:differential-window}. 

Since $g\in M^1(\R)$, failure of the frame
property and \Cref{prop:pregramian} yield $x\in\R$ and a nonzero
$\bm{c}=(c_k)_{k\in\Z}\in\ell^\infty(\Z)$ satisfying
\eqref{eq:bounded-coefficient-relation}. Set $x_*=bx$ and define $F$ by
\eqref{eq:gaussian-entire-function}.

For $z=t+is$, we have $|\phi_\beta(t+is-k)|=\myexp{\beta s^2}\myexp{-\beta(t-k)^2}$.
Using $|R_r(z)|\le C'_{r,\beta}(1+|z|)^r$ and
$1+|t-k+is|\le(1+|s|)(1+|t-k|)$, we obtain
\[
  |\phi_\beta^{(r)}(t+is-k)|
  \le C'_{r,\beta}(1+|s|)^r(1+|t-k|)^r
          \myexp{\beta s^2}\myexp{-\beta(t-k)^2}.
\]
Since the $1$-periodic continuous function
$\sum_{k\in\Z}(1+|t-k|)^r \myexp{-\beta(t-k)^2}$ is bounded on $\R$, the
series for $F$ and its termwise derivatives converge absolutely and
locally uniformly on $\C$, so $F$ is entire and satisfies \eqref{eq:F-growth}.

To verify that $F\not\equiv0$, write
\[
  \myexp{\beta z^2}F(z)
  =\sum_{k\in\Z}c_k \myexp{-\beta k^2}\myexp{2\beta kz}.
\]
If $F\equiv0$, evaluation at $z=i\theta/(2\beta)$ yields an identically
zero, absolutely convergent Fourier series
$\sum_{k\in\Z}c_k \myexp{-\beta k^2}\myexp{ik\theta}$, which implies $\bm{c}=\bm{0}$
by uniqueness of Fourier coefficients, contradicting $\bm{c}\ne\bm{0}$.

Since $Q(D)$ commutes with translations, termwise differentiation,
\eqref{eq:differential-window}, and \eqref{eq:bounded-coefficient-relation} give
\[
  Q(D)F(x_*+j\alpha)
  =\sum_{k\in\Z}c_k (Q(D)\phi_\beta)(x_*+j\alpha-k)
  =\sum_{k\in\Z}c_k g(x+aj-k/b)=0,
\]
proving \eqref{eq:lattice-differential-relation}. Finally,
\[
  \myexp{-\beta(z+i\pi/\beta-k)^2}
  =\myexp{\pi^2/\beta-2\pi iz}\myexp{2\pi ik}\myexp{-\beta(z-k)^2}
  =\myexp{\pi^2/\beta-2\pi iz}\myexp{-\beta(z-k)^2},
\]
thus multiplying by $c_k$ and summing 
over $k\in\Z$ gives \eqref{eq:F-automorphy}.
\end{proof}

\begin{example}[Hermite windows]
\label{ex:hermite}
Let $\phi(t)=\myexp{-\pi t^2}$. Using $D\phi=-2\pi t\phi$ and
$D^2\phi=(4\pi^2t^2-2\pi)\phi$, the first three Hermite functions in
\eqref{eq:hermite} can be written via constant-coefficient differential
operators:
\[
  h_0=\kappa_0\phi,\qquad
  h_1=-2\kappa_1D\phi,\qquad
  h_2=\kappa_2(16\pi^2t^2-4\pi)\phi=4\kappa_2(D^2+\pi)\phi.
\]
For the scaled window $h_n(t/b)$ in \Cref{prop:gaussian-reduction} with
$\beta=\pi/b^2$, the chain rule $D^r[\phi(t/b)]=b^{-r}\phi^{(r)}(t/b)$
gives the corresponding operators
\[
  Q_0(D)=\kappa_0,\qquad
  Q_1(D)=-2b\kappa_1D,\qquad
  Q_2(D)=4\kappa_2(b^2D^2+\pi).
\]
Thus \eqref{eq:lattice-differential-relation} reduces to $F'=0$ on the
shifted lattice $x_*+\alpha\Z$ for $h_1$, whereas for $h_2$ it becomes
the second-order relation $b^2F''+\pi F=0$.
\end{example}

\section{Wronskians and an order-\texorpdfstring{$n$}{n} differential relation}
\label{sec:wronskian}

For analytic functions $f_0,\ldots,f_{N-1}$, write
\begin{equation*}
 \Wr(f_0,\ldots,f_{N-1})(z)
       =\det\bigl(f_j^{(r)}(z)\bigr)_{r,j=0}^{N-1}.
\end{equation*}
The following elementary lemma extends the common-critical-point
argument in \cite[Lemma~3.5]{FaulhuberPetersen2026} to a differential
relation of order $n$.

\begin{lemma}[Wronskian multiplicity]
 \label{lem:multiplicity}
 Let $0\le n<N$, and let $f_0,\ldots,f_{N-1}$ be analytic near
 $\zeta\in\C$. Suppose there exist complex numbers $q_0,\ldots,q_n$, not all zero, such that
 \begin{equation}
  \sum_{r=0}^n q_r f_j^{(r)}(\zeta)=0,
               \qquad \text{for all } j \in \{0,\ldots,N-1\}.
  \label{eq:common-relation}
 \end{equation}
 Then
 \begin{equation}
  \frac{d^\ell}{dz^\ell}\Wr(f_0,\ldots,f_{N-1})(\zeta)=0,
                   \qquad \text{for } 0\le\ell<N-n.
  \label{eq:multiplicity-derivatives}
 \end{equation}
 In particular, if the Wronskian is not identically zero, its zero at
 $\zeta$ has multiplicity at least $N-n$.
\end{lemma}
\begin{proof}
 Define the row vectors
 \begin{equation*}
  v_r=\bigl(f_0^{(r)}(\zeta),\ldots,
                   f_{N-1}^{(r)}(\zeta)\bigr),\qquad r\ge0.
 \end{equation*}
 Condition \eqref{eq:common-relation} gives that the vectors 
 $v_0,\ldots,v_n$
 are linearly dependent. Fix $0\le\ell<N-n$. The determinant
 product rule gives
 \begin{equation*}
  \frac{d^\ell}{dz^\ell}\Wr(f_0,\ldots,f_{N-1})(\zeta)
  =\sum_{\substack{d_0,\ldots,d_{N-1}\in\Nzero\\d_0+\cdots+d_{N-1}=\ell}}
    \frac{\ell!}{d_0!\cdots d_{N-1}!}
    \det\bigl(v_{r+d_r}\bigr)_{r=0}^{N-1}.
 \end{equation*}
 Here $\det\bigl(v_{r+d_r}\bigr)_{r=0}^{N-1}$ denotes the determinant
 of the $N\times N$ matrix with row vectors $v_k$ indexed by $k=d_0, 1+d_1, \ldots, N-1+d_{N-1}$. A repeated entry in the list $d_0, 1+d_1, \ldots, N-1+d_{N-1}$ makes the corresponding determinant $\det\bigl(v_{r+d_r}\bigr)_{r=0}^{N-1}$ zero.
 Otherwise the indices $r+d_r$  for $r=0,\ldots,N-1$ are $N$ distinct nonnegative integers
 with sum
 \begin{equation}
  \sum_{r=0}^{N-1}(r+d_r)= \sum_{r=0}^{N-1} r + \sum_{r=0}^{N-1} d_r =
  \frac{N(N-1)}2+\ell.
  \label{eq:index-sum}
 \end{equation}
 The list $r+d_r$, $r=0,\ldots,N-1$, must contain every integer $0,\ldots,n$. Indeed, if $s\le n$
 were missing, the smallest possible sum of $N$ distinct nonnegative
 integers avoiding $s$ would be the sum of $0,1,\dots, s-1,s+1,\dots,N$, which is
 \begin{equation*}
  \sum_{r=0}^N r-s
   =\frac{N(N-1)}2+N-s
   \ge\frac{N(N-1)}2+N-n,
 \end{equation*}
 contradicting \eqref{eq:index-sum}. Each remaining determinant
 therefore contains all the linearly dependent rows
 $v_0,\ldots,v_n$ and vanishes. This proves
 \eqref{eq:multiplicity-derivatives}.
\end{proof}

\begin{remark}
 The local multiplicity $N-n$ in \eqref{eq:multiplicity-derivatives} is sharp. Take the $N$ monomials $z^j$ whose
 exponents are $j=0,\ldots,N$ with $n$ omitted. Their $n$th derivatives
 vanish at zero, hence $q_n=1$ and $q_r=0$ for $r\ne n$ satisfy \eqref{eq:common-relation}. Their Wronskian is a nonzero constant times
 $z^{N-n}$ (see \cite[Lemma 1]{BostanWronskians2010}), thus the zero at $0$ has multiplicity exactly $N-n$.
\end{remark}

Let $F$ be the function in \Cref{prop:gaussian-reduction}. For an integer $N\ge1$,
define
\begin{equation*}
 F_j(z)=F(z-j\alpha),\qquad j=0,\ldots,N-1,\qquad
 W_N=\Wr(F_0,\ldots,F_{N-1}).
\end{equation*}
The imaginary quasi-periodicity of $F$ both establishes independence
of these translates and propagates the zeros of their Wronskian.

\begin{lemma}
 \label{lem:WN}
Assume the setup of \Cref{prop:gaussian-reduction}.
Suppose that $N \in \N$ is such that
 \begin{equation}
  1,\myexp{2\pi i\alpha},\ldots,\myexp{2\pi i(N-1)\alpha}
       \quad\text{are pairwise distinct}.
  \label{eq:distinct-characters}
 \end{equation}
 Then $W_N\not\equiv0$, and
 \begin{equation}
  W_N\Bigl(z+\frac{i\pi}{\beta}\Bigr)
   =\myexp{N\pi^2/\beta-2\pi iNz+\pi i\alpha N(N-1)}W_N(z).
  \label{eq:WN-automorphy}
 \end{equation}
 Moreover,
 \begin{equation}
  \abs{W_N(t+is)}
  \le C_N(1+\abs s)^{N(N-1)/2}\myexp{N\beta s^2},
                           \qquad t,s\in\R.
  \label{eq:WN-growth}
 \end{equation}
 If $N>n$, the differential relation
 \eqref{eq:lattice-differential-relation} further implies
 \begin{equation}
  \ord_{\zeta}W_N\ge N-n,
   \qquad \zeta\in x_*+\alpha\Z+\frac{i\pi}{\beta}\Z.
  \label{eq:WN-lattice}
 \end{equation}
\end{lemma}
\begin{proof}
  Put $\rho(z)=\myexp{\pi^2/\beta-2\pi iz}$ and $\chi_j=\myexp{2\pi ij\alpha}$ for $j=0,\ldots,N-1$. Then \eqref{eq:F-automorphy} gives
  \begin{equation}
    \label{eq:F-automorphy-shifted}
  F_j\Bigl(z+\frac{i\pi}{\beta}\Bigr)
  =F\Bigl(z+\frac{i\pi}{\beta}-j\alpha\Bigr)
  =\myexp{\pi^2/\beta-2\pi i(z-j\alpha)}F(z-j\alpha)
  =\rho(z)\chi_j F_j(z),
  \end{equation}  
  Define an operator on entire functions by
  \begin{equation}
  (V_\beta f)(z)=\rho(z)^{-1} f\Bigl(z+\frac{i\pi}{\beta}\Bigr).
  \end{equation}
  With this definition, the quasi-periodicity \eqref{eq:F-automorphy-shifted} becomes
  $V_\beta F_j=\chi_j F_j$. Thus,
  the functions $F_j$ are nonzero eigenvectors with distinct
  eigenvalues, so they are linearly independent. Recall the classical result that a finite family of linearly independent (real or complex) analytic functions has a nonzero Wronskian, see, for example, \cite{BostanWronskians2010}. Hence, $W_N\not\equiv0$.

  To prove the quasi-periodicity of $W_N$, \eqref{eq:WN-automorphy},
  we follow \cite[Lemma~4.1]{FaulhuberPetersen2026}.
  By Leibniz' rule,
  \[
    (\rho F_j)^{(r)}=\sum_{k=0}^r\binom{r}{k}\rho^{(r-k)}F_j^{(k)},
  \]
  so the Wronskian matrix of $(\rho F_0,\ldots,\rho F_{N-1})$ is obtained
  from that of $(F_0,\ldots,F_{N-1})$ by left-multiplication by a lower
  triangular matrix with diagonal entries $\rho$. Thus
  \[
    \Wr(\rho F_0,\ldots,\rho F_{N-1})=\rho^N\Wr(F_0,\ldots,F_{N-1}).
  \]
  Factoring out the column scalars $\chi_j$, whose product is
  $\prod_{j=0}^{N-1}\chi_j=\myexp{\pi i\alpha N(N-1)}$, gives \eqref{eq:WN-automorphy}.

  Since \eqref{eq:F-growth} is uniform in $t\in\R$, each entry in row $r$ of
  $W_N(t+is)$ is bounded by $C_{r,\beta}\|\bm{c}\|_{\ell^\infty}(1+|s|)^r \myexp{\beta s^2}$.
  Expanding the determinant therefore gives \eqref{eq:WN-growth}, as each
  product takes one factor from every row $r=0,\ldots,N-1$ and
  $\sum_{r=0}^{N-1}r=N(N-1)/2$.

  Finally, fix $\zeta=x_*+k\alpha$ for an integer $k$. The operator $Q(D)$ commutes
  with real translation, and therefore
  \begin{equation*}
  \sum_{r=0}^n q_r F_j^{(r)}(\zeta)
    =Q(D)F\bigl(x_*+(k-j)\alpha\bigr)=0
        \qquad \text{for all } j\in\{0,\ldots,N-1\}.
  \end{equation*}
By \Cref{lem:multiplicity}, the Wronskian $W_N$ has a zero of multiplicity at least $N-n$ at
every real lattice point in $x_*+\alpha\Z$. The multiplier in
  \eqref{eq:WN-automorphy} is entire and nowhere zero, so translation
  by $i\pi/\beta$ preserves zero multiplicity. Iterating in both
  directions yields \eqref{eq:WN-lattice}.
\end{proof}

\section[Proofs of Theorems 1 and 2]{Proofs of \Cref{thm:irrational,thm:rational}}
\label{sec:proofs}
We compare the number of zeros in \eqref{eq:WN-lattice} with the
Gaussian growth bound \eqref{eq:WN-growth}. The required estimate is
a direct application of Jensen's formula.

\begin{lemma}
 \label{lem:zero-count}
 Let $H$ be a nonzero entire function. Suppose that for constants
 $L,\beta,C>0$ and $M\ge0$,
 \begin{equation}
  \abs{H(t+is)}\le C(1+\abs s)^M \myexp{L\beta s^2},
                          \qquad t,s\in\R.
  \label{eq:general-growth}
 \end{equation}
 If, for some $x_*\in\R$, $\alpha>0$, and integer $m\ge1$, all
 points of $x_*+\alpha\Z+\frac{i\pi}{\beta}\Z$ are zeros of $H$ of
 multiplicity at least $m$, then
 \begin{equation}
  \frac m\alpha\le L.
  \label{eq:zero-density-inequality}
 \end{equation}
\end{lemma}
\begin{proof}
 Choose $z_0\in\C$ with $H(z_0)\ne0$, and let $n_H(r)$ count the
 zeros, with multiplicity, of $H$ in the disk $\abs{z-z_0}<r$.
 Counting cells of the rectangular lattice gives
 \begin{equation*}
  \#\bigl((x_*+\alpha\Z+\frac{i\pi}{\beta}\Z)\cap B(z_0,r)\bigr)
       =\frac{\beta r^2}{\alpha}+O(r+1).
 \end{equation*}
 Indeed, cells fully enclosed by the disk and those intersecting its boundary 
 differ only by a boundary layer of bounded width and area $O(r)$.
 
 Since each shifted lattice point in $x_*+\alpha\Z+\frac{i\pi}{\beta}\Z$ is a zero of $H$ of multiplicity at least $m$, 
 we have for all $r>0$,
 \begin{align}
  n_H(r)&\ge\frac{m\beta}{\alpha}r^2-O(r+1), \nonumber 
  \intertext{which, after integration for $r \ge 1$, yields as $R\to\infty$}
  \int_1^R\frac{n_H(r)}r\,dr
     &\ge\frac{m\beta}{2\alpha}R^2-O(R).
     \label{eq:jensen-lower}
 \end{align}
 Restricting the integral to $r \ge 1$ omits only a bounded interval near zero that contributes a finite additive constant.

 For radii $R$ whose boundary contains no zero, Jensen's formula and
 \eqref{eq:general-growth} give
 \begin{align}
  \int_0^R\frac{n_H(r)}r\,dr
   &=\frac1{2\pi}\int_0^{2\pi}
          \log\abs{H(z_0+R\myexp{i\theta})}\,d\theta-\log\abs{H(z_0)} \nonumber \\
   &\le\frac{L\beta}{2}R^2+O(\log R),
   \label{eq:jensen-upper}
 \end{align}
 where we used
 \begin{equation*}
  \frac1{2\pi}\int_0^{2\pi}
      (\operatorname{Im}z_0+R\sin\theta)^2\,d\theta
       =(\operatorname{Im}z_0)^2+\frac{R^2}{2}.
 \end{equation*}
 Comparing \eqref{eq:jensen-lower} and \eqref{eq:jensen-upper},
 dividing by $\beta R^2/2$, and letting $R$ tend to infinity
 through admissible radii (those for which the boundary contains no zero)
 proves \eqref{eq:zero-density-inequality}.
\end{proof}

\begin{proposition}
 \label{prop:obstruction}
 Let $g$ be a polynomial--Gaussian window
 as in \eqref{eq:window} with $\deg P=n$.
 Suppose $\G(g,a,b)$ is not a frame, and let $\alpha=ab$.
 For every integer $N>n$ satisfying
 \eqref{eq:distinct-characters}, one has
 \begin{equation*}
  N-n\le N\alpha.
 \end{equation*}
\end{proposition}
\begin{proof}
 \Cref{prop:gaussian-reduction} supplies a nonzero entire function $F$
 of the form \eqref{eq:gaussian-entire-function} and the lattice relation
 \eqref{eq:lattice-differential-relation}. By \Cref{lem:WN}, the function $W_N$
 is nonzero, obeys \eqref{eq:WN-growth}, and has multiplicity at
 least $N-n$ at every point of
 $x_*+\alpha\Z+\frac{i\pi}{\beta}\Z$.
 Applying \Cref{lem:zero-count} with $H=W_N$, $L=N$, and $m=N-n$ completes the proof.
\end{proof}

\begin{proof}[Proof of \Cref{thm:irrational}]
 Suppose towards a contradiction that $\G(g,a,b)$ is not a frame for 
 irrational $\alpha=ab \in \itvoo{0,1}$. Since $\alpha$ is irrational, the exponential terms 
 \eqref{eq:distinct-characters} are distinct for every $N$.
 Thus \Cref{prop:obstruction} gives
 \begin{equation*}
  1-\frac nN\le\alpha\qquad\text{for every integer }N>n.
 \end{equation*}
 Letting $N\to\infty$ yields $1\le\alpha$, which contradicts the assumption $\alpha<1$ in \eqref{eq:irrational-condition}. Therefore $\G(g,a,b)$ is a frame.
\end{proof}

\begin{proof}[Proof of \Cref{thm:rational}]
 Suppose towards a contradiction that $\G(g,a,b)$ is not a frame
 for $ab=p/q$ with $\gcd(p,q)=1$ and write $\alpha=p/q$.
The inequality $q-p>n$ implies
 $q>n$. Since $p$ and $q$ are relatively prime, the numbers
$\myexp{2\pi ij\alpha}=\myexp{2\pi ijp/q}$, $j=0,\ldots,q-1$ are distinct. Taking $N=q$ in \Cref{prop:obstruction} gives $q-n\le q\alpha=p$, which contradicts $q-p>n$. Hence the system is a frame.
\end{proof}

\section[Declaration of use of AI tools]{Declaration of use of generative AI and AI-assisted technologies}

During the preparation of this manuscript, the author used ChatGPT
(OpenAI) as an interactive tool for exploring mathematical
arguments, checking intermediate calculations, developing and organizing
proofs, and improving the presentation. All arguments, results, and references included in the manuscript were verified by the author, who takes full responsibility for the final content.


\end{document}